\documentclass[11pt]{article}
\usepackage{amssymb}
\usepackage[pagebackref=true,colorlinks,linkcolor=red,citecolor=blue]{hyperref}
\newtheorem{precor}{{\bf Corollary}}

\newenvironment{cor}{\begin{precor}{\hspace{-0.5
               em}{\bf.\ }}}{\end{precor}}
\newtheorem{prerem}{{\bf Remark}}

\newenvironment{rem}{\begin{prerem}{\hspace{-0.5
               em}{\bf.\ }}}{\end{prerem}}
\newtheorem{precon}{{\bf Conjecture}}

\newtheorem{predefin}{{\bf Definition}}

\newenvironment{defin}[1]{\begin{predefin}{\hspace{-0.5
                   em}{\bf.\ }}{\rm
#1}\hfill{$\spadesuit$}}{\end{predefin}}
\newtheorem{preexm}{{\bf Example}}

\newtheorem{preappl}{{\bf Application}}

\newtheorem{prelem}{{\bf Lemma}}

\newenvironment{lem}{\begin{prelem}{\hspace{-0.5
               em}{\bf.\ }}}{\end{prelem}}
\newtheorem{preproof}{{\bf Proof.\ }}

\newenvironment{proof}[1]{\begin{preproof}{\rm
               #1}\hfill{$\blacksquare$}}{\end{preproof}}
\newtheorem{presproof}{{\bf Sketch of Proof.\ }}

\newtheorem{prethm}{{\bf Theorem}}

\newenvironment{thm}{\begin{prethm}{\hspace{-0.5
               em}{\bf.\ }}}{\end{prethm}}
\newtheorem{prealphthm}{{\bf Theorem}}

\newenvironment{alphthm}{\begin{prealphthm}{\hspace{-0.5
               em}{\bf.\ }}}{\end{prealphthm}}
\newtheorem{prealphlem}{{\bf  Lemma}}

\newenvironment{alphlem}{\begin{prealphlem}{\hspace{-0.5
               em}{\bf.\ }}}{\end{prealphlem}}
\newtheorem{prepro}{{\bf Proposition}}

\newtheorem{preprb}{{\bf Problem}}

\newtheorem{prealphprb}{{\bf Problem}}

\newtheorem{prequ}{{\bf Question}}

\def\conct[#1,#2]{\mbox {${#1} \leftrightarrow {#2}$}}
\def\dconct[#1,#2]{\mbox {${#1} \rightarrow {#2}$}}
\def\deg[#1,#2]{\mbox {$d_{_{#1}}(#2)$}}
\def\mindeg[#1]{\mbox {$\delta_{_{#1}}$}}
\def\maxdeg[#1]{\mbox {$\Delta_{_{#1}}$}}
\def\outdeg[#1,#2]{\mbox {$d_{_{#1}}^{^+}(#2)$}}
\def\minoutdeg[#1]{\mbox {$\delta_{_{#1}}^{^+}$}}
\def\maxoutdeg[#1]{\mbox {$\Delta_{_{#1}}^{^+}$}}
\def\indeg[#1,#2]{\mbox {$d_{_{#1}}^{^-}(#2)$}}
\def\minindeg[#1]{\mbox {$\delta_{_{#1}}^{^-}$}}
\def\maxindeg[#1]{\mbox {$\Delta_{_{#1}}^{^-}$}}
\def\isdef{\mbox {$\ \stackrel{\rm def}{=} \ $}}
\def\dre[#1,#2,#3]{\mbox {${\cal E}_{_{#3}}(#1,#2)$}}
\def\pdre[#1,#2,#3]{\mbox {${\cal P}_{_{#3}}(#1,#2)$}}
\def\var[#1,#2]{\mbox {${\rm Var}_{_{#1}}(#2)$}}
\def\ls[#1]{\mbox {$\xi^{^{#1}}$}}
\def\hom[#1,#2]{\mbox {${\rm Hom}({#1},{#2})$}}
\def\onvhom[#1,#2]{\mbox {${\rm Hom^{v}}(#1,#2)$}}
\def\onehom[#1,#2]{\mbox {${\rm Hom^{e}}(#1,#2)$}}
\def\core[#1]{\mbox {$#1^{^{\bullet}}$}}
\def\cay[#1,#2]{\mbox {${\rm Cay}({#1},{#2})$}}
\def\cays[#1,#2]{\mbox {${\rm Cay_{s}}({#1},{#2})$}}
\def\dirc[#1]{\mbox {$\stackrel{\rightarrow}{C}_{_{#1}}$}}
\def\cycl[#1]{\mbox {${\bf Z}_{_{#1}}$}}

\begin{document}
\footnotetext[1]{$\ast$This research was in part supported by
a grant from IPM (No. 89050112).}
\begin{center}
{\Large \bf Some New Bounds For Cover-Free Families\\ Through  Biclique Cover}\\
\vspace*{0.5cm}
{\bf Hossein Hajiabolhassan$^\ast$ and Farokhlagha Moazami$^\dag$}\\
{\it $^\ast$Department of Mathematical Sciences}\\
{\it Shahid Beheshti University, G.C.}\\
{\it P.O. Box {\rm 1983963113}, Tehran, Iran}\\
{\it $^\ast$School of Mathematics\\
Institute for Research in Fundamental Sciences {\rm (}IPM{\rm )}}\\
{\it P.O. Box {\rm 193955746}, Tehran, Iran}\\
{\tt hhaji@sbu.ac.ir}\\
{\it $^\dag$Department of Mathematics} \\
{\it Alzahra University}\\
{\it P.O. Box {\rm 1993891176}, Tehran, Iran}\\
{\tt f.moazami@alzahra.ac.ir}\\ \ \\
\end{center}
\begin{abstract}
An $(r,w;d)$ cover-free family $(CFF)$ is a family of subsets of
a finite set such that the intersection of any $r$ members of the
family contains at least $d$ elements that are~not in the union of
any other $w$ members. The minimum number of elements for which
there exists an $(r,w;d)-CFF$ with $t$ blocks is denoted by
$N((r,w;d),t)$.

In this paper, we show that the value of $N((r,w;d),t)$ is equal
to the $d$-biclique covering number of the bipartite graph
$I_t(r,w)$ whose vertices are all $w$- and $r$-subsets of a
$t$-element set, where a $w$-subset is adjacent to an $r$-subset
if their intersection is empty. Next, we introduce some new bounds for
$N((r,w;d),t)$. For instance, we show that for $r\geq w$ and $r\geq 2$
$$ N((r,w;1),t) \geq c{{r+w\choose w+1}+{r+w-1 \choose w+1}+ 3 {r+w-4 \choose w-2} \over \log r} \log (t-w+1),$$
where $c$ is a constant satisfies the well-known bound $N((r,1;1),t)\geq c\frac{r^2}{\log r}\log t$.
Also, we determine the exact value of
$N((r,w;d),t)$ for some values of $d$. Finally, we show that $N((1,1;d),4d-1)=4d-1$ whenever
there exists a Hadamard matrix of order $4d$.
\begin{itemize}
\item[]{{\footnotesize {\bf Key words:}\ cover-free family, biclique cover,
fractional biclique cover, weakly cross-intersecting set-pairs.}}
\item[]{ {\footnotesize {\bf Subject classification:} 05B40.}}
\end{itemize}
\end{abstract}
\section{Introduction}
A family of sets is called an $(r,w)$-cover-free family if no
intersection of $r$ sets of the family are covered by a union of
any other $w$ sets of the family. Cover-free family was first
introduced by Kautz and Singleton \cite{first} to investigate the
properties of the non-random binary superimposed codes. In
$1985$, Erd\"{o}s, Frankl, and F\"{u}redi \cite{erdos} studied
this concept as a generalization of Sperner system. In $1988$,
Mitchell and Piper \cite{kps} defined the concept of key
distribution pattern which is in fact a generalized type of
cover-free family. Others have used this concept in cryptography,
for example, group key predistribution, frameproof codes, broad
cast anti-jamming, and so on, see \cite{colborn}. Cover-free family
has been studied extensively throughout the literature
due to both its interesting structure and the central role it plays in several respects, see
\cite{blackburn, engel, erdos, first, cff2, wei}.

In this paper, we discuss aspects relevant to cover-free families.
In Section $1$, we set up notation and terminology. Section $2$
is  devoted to study the connection between cover-free families
and biclique cover. In Section $3$, we presents several new lower bounds for
$N((r,w;d),t)$. Finally, Section $4$ concerns
the fractional version of biclique cover and we determine the
exact value of $N((r,w;d),t)$ for some values of $d$.
Finally, we show that if there exists a Hadamard matrix of order $4d$, then $N((1,1;d),4d-1)=4d-1$.

Throughout this paper, we only consider finite simple graphs. For
a graph $G$, let $V(G)$ and $E(G)$ denote its vertex and edge
sets, respectively. By a {\it biclique} we mean a bipartite graph
with vertex set $(X ,Y)$ such that every vertex in $X$ is
adjacent to every vertex in $Y$. Note that every empty graph is a
biclique. A {\it biclique  cover} of a graph $G$ is a collection
of bicliques of $G$ such that each edge of $G$ is in at least one
of the bicliques. The number of bicliques in a minimum biclique
cover of $G$ is called the { \em biclique  covering  number} of
$G$ and denoted by $bc(G)$. This measure of graphs is studied in
the literature \cite{ bc3, bc1, bc2}.

In this paper, we also need a generalization of biclique cover as
follows.
\begin{defin}{A  $d$-{\it biclique cover} of a graph $G$ is a collection of bicliques
  of $G$ such that each edge of
$G$ is in at least $d$ of the bicliques. The number of bicliques
in a minimum $d$-biclique cover of $G$ is called the {\em
$d$-biclique covering  number} of $G$ and denoted by $bc_d(G)$.}
\end{defin}
As usual, we denote by $[t]$ the set $\{1, 2,
\ldots, t\}$. In this
paper, by $A^c$ we mean the complement of the set $A$.
 For $0 < w \leq r \leq t$, the {\it subset
 graph} $S_t(w,r)$ is a bipartite graph whose vertices are all
$w$- and $r$-subsets of a $t$-element set, where a $w$-subset is
adjacent to an $r$-subset if and only if one subset is contained
in the other. Some properties of this family of graphs have been
studied by several researchers, see \cite{ subset2, subset1,
fractional}. In this paper, we consider an isomorphic version of
this graph and name it {\em bi-intersection graph}.
 \begin{defin}{ For $0 < w \leq r \leq t$, the {\it bi-intersection graph} $I_t(r,w)$ is a
 bipartite graph whose vertices are all
$w$- and $r$-subsets of a $t$-element set, where a $w$-subset is
adjacent to an $r$-subset if and only if their intersection is
empty.}
  \end{defin}
 It is~not difficult to see that the bi-intersection graph  $I_t(r,w)$ is isomorphic to
 $S_t(r,t-w)$. A {\em set  system} is an ordered pair $(X,{\cal B})$, where $X$
is a set of elements and ${\cal B}$ is a family of subsets
(called block) of $X$. A set system can be described by an
incidence matrix. Let $(X, {\cal B})$ be a set system, where $X=\{
x_1,x_2, \ldots, x_v\}$ and ${\cal B}=\{ B_1,B_2, \ldots, B_b\}$.
The incidence matrix of $(X, {\cal B})$ is the $b\times v$ matrix
$A=(a_{ij})$, where $$ a_{ij}=\left \{
\begin{array}{cc}
1  & if \,\,\,\ x_j \in B_i\\
0 & if \,\,\,\ x_j \notin B_i\\
\end{array} \right .
$$
The next definition is a formal definition of cover-free family.

\begin{defin}{
Let $n$, $t$, $r$, and $w$ be positive integers. A set  system
$(X,{\cal B})$, where $|X|=n$ and ${\cal B}=\{B_1, \ldots,B_t\}$
is called   an $(r,w)-CFF(n,t)$ if for any two sets of indices
$L,M \subseteq[t]$ such that $L \cap M = \varnothing $, $ |L| = r
$, and $|M|=w$, we have
$$\displaystyle \bigcap_{l \in L} B_l \nsubseteq  \displaystyle \bigcup_{m \in M}B_m.$$}
\end{defin}

Stinson and Wei \cite{cff1} generalized the definition of cover-free family as
follows.
\begin{defin}{Let $d$, $n$, $t$, $r$, and $w$ be positive integers. A set system  $(X,{\cal B})$,
where $|X|=n$ and ${\cal B}=\{B_1, \ldots,B_t\}$ is called   an $(r,w;d)-CFF(n,t)$
if for any two sets of indices $L,M \subseteq[t]$ such that $L
\cap M = \varnothing $, $ |L| = r $, and $|M|=w$, we have
$$
|(\displaystyle \bigcap_{l \in L}B_l)\setminus
(\displaystyle\bigcup_{m \in M}B_m)| \geq d .$$}
\end{defin}

Let $N((r,w;d),t)$ denote the minimum number of elements in any
$(r,w;d)-CFF$ having $t$ blocks. For convenience, we use the
notation $N((r,w),t)$ instead of $N((r,w;1),t)$. Obviously, we
have $N((r,w;d),t)=N((w,r;d),t)$. Hence, unless otherwise stated
we assume that $w\leq r$.
\section{Biclique Cover}
In this section, we show that the existence of a cover-free
family can result from the existence of biclique cover of
bi-intersection graph and vice versa. Our viewpoint sheds some new light on
cover-free family. Using this observation, we
introduce several new bounds.

 \begin{thm}\label{equi} Let $r$, $w$, and $t$ be  positive integers, where $t \geq
 r+w$. It holds that
 $$N((r,w),t)=bc(I_t(r,w)).$$
 \end{thm}
\begin{proof}{
First, consider an optimal $(r,w)-CFF(n,t)$, i.e.,
$n=N((r,w),t)$, with incidence matrix $A=(a_{ij})$. Assign to the
$j$th column of $A$,  the set $A_j$ as follows
$$A_j\isdef \{i|\  1\leq i \leq t, a_{ij}=1\}.$$

 Now, for any $1\leq j \leq n$, construct a bipartite graph $G_j$ with vertex set $(X_j,Y_j)$, where the
vertices of $X_j$ are all $r$-subsets of $A_j$ and the vertices
of $Y_j$ are all $w$-subsets of $A_j^c$, i.e.,
$$ X_j=\{ U | \,\,\ U\subseteq A_j, \,\ |U|=r \}\quad {\rm and} \quad Y_j=\{V |\,\,\ V \subseteq A_j^c, \,\ |V|=w\}.$$
Also, an $r$-subset is adjacent to a $w$-subset if their
intersection is empty. One can see that $G_j$, for $1\leq j\leq n$,
is a biclique. Let $UV$ be an arbitrary edge of $I_t(r,w)$,
where $U \cap V = \varnothing$, $|U|=r$ and $|V|=w$. In view of
definition of CFF and since $A$ is the incidence matrix of the
$CFF$, there is a column of $A$, say $j$, where $a_{ij}=1$ if
$i\in U$ and $a_{ij}=0$ if $i\in V$. Clearly, $U \in X_j$, $V \in
Y_j$, and $UV \in G_j$. Hence, $\{G_1,G_2, \ldots, G_n\}$ is a
biclique cover of $I_t(r,w)$. So $bc(I_t(r,w)) \leq N((r,w),t)$.

Conversely, assume that $G_1, \ldots, G_l$ constitute a biclique
cover of $I_t(r,w)$, where $l=bc(I_t(r,w))$ and $G_i$ has as its
vertex set $(X_i, Y_i)$. Let $A_i$ be the union of sets that lie
in $X_i$. Consider the indicator vector of the set $A_i$, for $i
= 1, \ldots, l$, and construct the matrix $A$ whose $i$th column
is the indicator vector of the set $A_i$.  We claim that $A$ is
the incidence matrix of an $(r,w)-CFF(l,t)$. To see this, let $U$
and $V$ be two arbitrary disjoint sets of $[t]$, where $|U|=r$ and
$|V|=w$. Thus, $UV$ is an edge of the graph $I_t(r,w)$. Hence,
there exists a biclique $G_j$, where $U \in X_j$ and  $V \in
Y_j$. Now, in view of definition of $A_j$, one can see that all
entries corresponding to the elements of $U$ and $V$ in the $j$th
column are $1$ and $0$, respectively. So $ N((r,w),t) \leq
bc(I_t(r,w))$. This completes the proof.}
\end{proof}
By the same argument we obtain the following corollary.

\begin{cor} Let $r$, $w$, $d$, and $t$ be  positive integers, where $t\geq r+w$. It holds that
$$N((r,w;d),t)=bc_d(I_t(r,w)).$$
\end{cor}
A {\em weakly separating system} on $[t]$ is a collection $\{(X_1,
Y_1), \ldots , (X_n, Y_n)\}$ of disjoint pairs of subsets of $[t]$
such that for every $i, j \in [t]$ with $i \neq j$ there is a $k$
with either $i \in X_k$ and $j \in Y_k$ or $i \in Y_k$ and $j \in
X_k$. Similarly, a {\em strongly separating system} on $[t]$ is a
collection $\{(X_1, Y_1),\ldots , (X_n, Y_n)\}$ of disjoint pairs
of subsets of $[t]$ such that for every ordered pair $(i,j)$
with $1\leq i, j \leq t$ and $i \neq j$, there is a $k\in [n]$ with
$i \in X_k$ and $j \in Y_k$. The study of separating systems was
started by R\'{e}nyi \cite{ss} in 1961. Other researchers have
studied the properties of separating system in the literature, see
\cite{bllobas1, bollobas2, spencer, css}. One can construct a
$(1,1)-CFF(n,t)$ from a strongly separating system on $[t]$ of
size $n$ and vice versa (see the proof of Theorem \ref{equi}). So
if we denote by ${\cal R}(t)$, the minimum size of a strongly
separating system, we have $N((1,1),t)={\cal R} (t)$. Let
$\{(X_1,Y_1), \ldots, (X_n,Y_n)\}$ be a weakly separating system.
The complete bipartite graphs with vertex classes $X_i$ and $Y_i$
cover the edges of the complete graph $K_t$ with vertex set
$[t]$. Also, if the family $\{G_1, \ldots, G_n\}$ is a biclique
cover of $K_t$, where $G_i$ has as its vertex set $(X_i,Y_i)$,
then $\{(X_1,Y_1), \ldots, (X_n,Y_n)\}$ is a weakly separating
system. So if we denote by $s(t)$, the size of minimum weakly
separating system, then we have $s(t)=bc(K_t)$. Also, in
\cite{bc1}, it was proved that ${\cal R}(t) = bc(K_{t,t}^-)$,
where $K_{t,t}^-$ is the complete bipartite graph $K_{t,t}$ with
a perfect matching removed. The exact value of ${\cal R}(t)$ was
determined by Spencer \cite{spencer}.
\begin{alphthm}
\label{sss} {\rm \cite{spencer}} If $C=\min\{ c \,\ | \,\ {c
\choose \lfloor \frac{c}{2}\rfloor}\geq t \}$, then $C={\cal
R}(t)$.
\end{alphthm}

Theorem \ref{sss} implies
$${\cal R}(t)= \log_2 t +\frac{1}{2}\log_2\log_2 t + O(1).$$
It is simple to see that $bc(G)\geq m(G)$, where $m(G)$ is the
maximum size of induced matchings of $G$. Let ${\cal F}=
\{(A_i,B_i)\}_{i=1}^h$ be a family of pairs of subsets of an
arbitrary set. The family ${\cal F}$ is called an $(r,w)$-system
if for all $1\leq i \leq  h$, $|A_i|=r$, $|B_i|=w$, $A_i\cap
B_i=\varnothing$, and for all distinct $i,j$ with $1\leq i,j \leq
h$, $A_i\cap B_j\neq\varnothing$. Bollob\'{a}s \cite{bollobas3}
proved that the maximum size of an $(r,w)$-system is equal to
${r+w \choose r}$. Obviously, $m(I_t(r,w))$ is the maximum size of
an $(r,w)$-system, so $N((r,w),t)\geq {r+w \choose r}$.
A $covering$ of a graph $G$ is a subset $K$ of $V(G)$ such that
every edge of $G$ has at least one end in $K$. The number of
vertices in a minimum covering of $G$ is called the {\it
covering  number} of $G$ and denoted by $\beta(G)$. It is~not
difficult to see that the biclique covering number of a graph $G$
without $C_4$ as a subgraph is equal to the covering number of
$G$. So we have the following corollary.
\begin{cor} For any positive integers $r$, $w$, and $t$, where $t=r+w+1$ or $t=r+w$, we have
\begin{center}
 $N((r,w),t)= \min\{{ t \choose r}, {t \choose w}\}= {t \choose w}.$
\end{center}
\end{cor}
\begin{proof}{It is easy  to see that the graph $I_t(r,w)$, when $t=r+w+1$ or $t=r+w$,
 does~not contain $C_4$ as a subgraph. On the other
hand, it is  well-known  that for every bipartite graph, the
covering number is equal to the maximum size of matchings. Easily,
using Hall's Theorem, the maximum number of matching in this
graph is equal to $\min\{{t \choose r}, {t \choose w}\}$.}
\end{proof}
We should mention that it is known \cite{engel} that $N((r,w),t)={t \choose w}$ whenever
$t\leq r+w+{r\over w}$. As an interesting application of cover-free family, one can consider key predistribution scheme (KPS).
The specification structure of a KPS is the family of all disjoint pairs $(P, F)$ of subsets of
the set of users $U$ such that every user in $P$ must be able to compute a common key of $P$ that
will remain unknown to the coalition $F$. The above corollary gives the exact
value of the minimum number of the keys in a KPS,
constructed by a cover-free family, with $r+w+1$ users.
\section{Bounds}
In this section, we introduce several bounds for  $N((r,w;d),t)$.
Engel \cite{engel}, using the fractional matching and fractional
cover of ordered interval hypergraph, obtained the following bounds
\begin{alphthm}\label{engel}{\rm \cite{engel}} For any positive integers $r$,
$w$, and $t$, where $r\geq w$ and $t\geq r+w$, we have
$$ N((r,w),t)\geq {r+w-1 \choose r} {\cal R}(t-r-w+2).$$
\end{alphthm}
\begin{alphthm}{\rm \cite{engel}} For any $\epsilon >0$, it holds that
$$ N((r,w),t_{\epsilon})\geq (1-\epsilon)\frac{(w+r-2)^{w+r-2}}{(w-1)^{w-1}(r-1)^{r-1}}{\cal R}(t_{\epsilon}-r-w+2),$$
for all sufficiently large $t_{\epsilon}$.
\end{alphthm}
Here is the best known lower bound for $N((r,1),t)$.
\begin{alphthm}\label{rcff}{\rm\cite{rcff1, rcff2, rcff3}} Let $r \geq 2$ and $t \geq r+1$ be  positive
integers. It holds that
$$N((r,1),t)\geq C_{r,t}\frac{r^2}{\log r}\log t,$$
where $\displaystyle \lim_{r+t\rightarrow \infty}C_{_{r,t}}=c$ for some constant $c$.
\end{alphthm}
 Several proofs have been presented for the preceding theorem. In \cite{rcff1, rcff2,
 rcff3}, it was shown that $c$ is approximately $\frac{1}{2}$,
$\frac{1}{4}$, and  $\frac{1}{8}$, respectively.

\begin{alphlem}\label{stinlem}{\rm \cite{cff2}} For any positive integers $r$,
$w$, and $t$, where $t\geq r+w$, we have
$$N((r,w),t)\geq N((r,w-1),t-1) + N((r-1,w),t-1).$$
\end{alphlem}

Stinson, Wei, and Zhu \cite{cff2}, using Lemma~\ref{stinlem}
and Theorem~\ref{rcff}, improved the bounds of Engel in some cases and obtained
the following bounds.
\begin{alphthm}\label{stinthm} {\rm \cite{cff2}}
For any positive integers  $r$, $w$, and $t$, where $t\geq r+w$, we have
$$ N((r,w),t) \geq 2c {{w+r \choose r}\over \log(w+r)} \log
t,$$
where $c$ is a constant satisfies Theorem~\ref{rcff}.
\end{alphthm}
\begin{alphthm}\label{stinthm2}{\rm \cite{cff2}} For any positive integers $r,w\geq 1$, there exists an integer $t_{r, w}$ such
that for all $t>t_{r, w}$
$$ N((r,w),t) \geq 0.7c (r+w) {{w+r \choose r}\over \log({w+r\choose r})} \log
t,$$
where $c$ is a constant satisfies Theorem~\ref{rcff}.
\end{alphthm}

In \cite{wei2}, it was shown $t_{r, w}\leq \max \{\lfloor {r+w+1 \over 2}\rfloor^2,5\}$.
Here we introduce some new lower bounds for $N((r,w;d),t)$ which improve Theorem~\ref{engel} and also
we present a lower bound (Theorem~\ref{estimate}) which can be considered as an improvement of Theorems~\ref{stinthm} and \ref{stinthm2} whenever $w$ is sufficiently small relative to $r$.
We first prove the following preliminary lemma which will be needed in the proof of Theorem~\ref{mainthm}.
\begin{lem}\label{abowe} Let $G$ be a graph and  $G_1,G_2, \ldots, G_k$ be some pairwise vertex disjoint subgraphs of $G$. Also, assume that for every four cycle $C_4$ of $G$ and $1\leq i\neq j\leq k$, we have $E(C_4)\cap E(G_i)=\varnothing$ or
$E(C_4)\cap E(G_j)= \varnothing$. Then
$$bc_d(G)\geq \sum_{i=1}^k bc_d(G_i).$$
\end{lem}
\begin{proof}{Let $\{ H_1, H_2, \ldots, H_l \}$ be an optimal $d$-biclique cover of $G$, i.e., $l=bc_d(G)$.
Let $H'_i$ be a subgraph of $G_1 \cup G_2 \cup \cdots \cup G_k$ induced by $H_i$. If
$H'_i$ is a non-empty graph, by the assumption, it is clear that
$H'_i$ is a biclique of exactly one of $G_i$'s. Now one can
see that $H'_j$'s cover all edges of $G_i$'s
at least $d$ times. So $bc_d(G)\geq \sum_{i=1}^k bc_d(G_i),$ as desired.}
\end{proof}
Before embarking on the proof of the next theorem, we need the following definition. The family ${\cal F}=\{(A_1,B_1), \ldots, (A_g,B_g)\}$ is called a {\it weakly cross-intersecting set-pairs} (resp. cross-intersecting set-pairs) of size $g$ on a ground set of cardinality $h$ whenever all $A_i$'s and $B_i$'s are subsets of an $h$-set and for every $i$, where $1\leq i \leq g$, $A_i$ and $B_i$ are disjoint subsets,
and furthermore, for every $i \neq j$, $(A_i \cap B_j)\cup (A_j \cap B_i)\not=\varnothing$ (resp. $(A_i \cap B_j)\not=\varnothing$ and $(A_j \cap B_i)\not=\varnothing$). This concept
is a variant of the generalization of $(r,w)$- weakly cross-intersecting set-pairs which was introduced first
by Tuza \cite{tuza}. The weakly cross-intersecting set-pairs ${\cal F}=\{(A_1,B_1), \ldots, (A_g,B_g)\}$ is
called an $(r,w)$-{\it weakly cross-intersecting set-pairs} whenever for every $1\leq i\leq g$, $|A_i|=r$
and $|B_i|=w$.
\begin{rem}{\rm It is worth noting that one can consider any weakly cross-intersecting set-pairs as a biclique covering. To see this, assume that $\{G_1, \ldots, G_l\}$ is a biclique cover of a graph $G$ and $G_i$ has $(X_i, Y_i)$, as its vertex set. This biclique cover is called an $(r,w)$-biclique cover whenever
each vertex of $G$ belongs to at most r sets in $\{X_1,X_2,\ldots,X_l\}$ and at most $w$ sets in  $\{Y_1,Y_2,\ldots,Y_l\}$.  Let  $\{(A_1, B_1), \ldots, (A_g, B_g) \}$ be a set-pairs on a ground set of size $h$. The dual set system $\{(S_1, T_1), \ldots, (S_h, T_h)\}$ is a set-pairs on a ground set of size $g$ and is defined by $S_i=\{ j : i \in A_j\}$ and  $T_i=\{ j : i \in B_j\}$, for $i=1, 2, \ldots,h $. It is shown in \cite{bollobas2} that a family is a
cross-intersecting set-pairs if and only if its dual is a strongly separating system. Similarly, one can see that $\{(A_1, B_1), \ldots, (A_g, B_g) \}$ is a weakly cross-intersecting set-pairs on a ground set of size $h$ such that for any $1\leq i\leq g$, $|A_i|\leq r$ and $|B_i|\leq w$ if and only if the dual of this set-pairs system, i.e.  $\{(S_1, T_1), \ldots, (S_h, T_h)\}$, is an $(r,w)$-biclique cover of size $h$ for the complete graph $K_g$.
}
\end{rem}
Hereafter, we adopt the convention that $N((r,0;d),t)=N((0,w;d),t)=1$.
\begin{thm}\label{mainthm} Suppose that $g, h, r, w$, and $t$ are positive integers. Also, let ${\cal F}=\{(A_1,B_1), \ldots, (A_g,B_g)\}$ be a weakly cross-intersecting
set-pairs on a ground set of size $h$ such that for any $1\leq i\leq g$, $|A_i|\leq r$ and $|B_i|\leq w$. If $t\geq \max\{h, r+w\}$, then
 $$N((r,w;d),t) \geq \sum_{i=1}^g N((r-|A_i|,w-|B_i|;d),t-|A_i|-|B_i|).$$
\end{thm}
\begin{proof}{Assume that ${\cal F}= \{(A_1,B_1), \ldots, (A_g,B_g)\}$ is a weakly cross-intersecting set-pairs. For every $1\leq k\leq g$, construct a bipartite graph $G_k$ with vertex set $(X_k,Y_k)$, where the vertices of $X_k$ are all $r$-subsets of $[t]$ which contain $A_k$ and their intersections with $B_k$ are empty.  Also, the vertices of $Y_k$ are all $w$-subsets of the set $[t]$ which contain $B_k$ and their
intersections with $A_k$ are empty, i.e.,
$$X_k = \{ U \, | \, U\subseteq [t], \, |U|=r, \,\,\ A_k \subseteq U, \,\,\ U \cap B_k = \varnothing \}$$
$$Y_k = \{ V \, | \, V\subseteq [t], \,  |V|=w, \,\,\ B_k \subseteq V, \,\,\ V \cap A_k = \varnothing \},$$
where a vertex $U \in X_k$ is adjacent to a vertex $V \in Y_k$ if
$U \cap V = \varnothing$. Clearly, if $|A_k|=r$ or $|B_k|=w$, then $G_k$ is isomorphic to a star graph. Otherwise,
one can check that every $G_k$ is isomorphic to $I_{t-|A_k|-|B_k|}(r-|A_k|,w-|B_k|)$. Since if we delete the elements of $A_k$ from the  vertices of $X_k$, every vertex is mapped  to an $(r-|A_k|)$-subset of the set  $[t]\setminus (A_k\cup B_k)$ and also if we
remove the elements of $B_k$ from the  vertices of $Y_k$, every vertex
is mapped  to a $(w-|B_k|)$-subset of the set $[t]\setminus (A_k\cup B_k)$.
Obviously, this mapping is an isomorphism between $G_k$ and
$I_{t-|A_k|-|B_k|}(r-|A_k|,w-|B_k|)$. On the other hand, since ${\cal F}$ is a weakly cross-intersecting set-pairs, $G_k$'s are pairwise vertex disjoint. Also, for any $1\leq i\neq j\leq k$, there is no four cycle $C_4$ of $I_{t}(r,w)$ such that $E(C_4)\cap E(G_i)\neq\varnothing$ and $E(C_4)\cap E(G_j)\neq \varnothing$. So, in view of Lemma~\ref{abowe}, one can see that $$bc_d(I_t(r,w)) \geq \sum _{k=1}^h bc_d(G_k).$$
Hence, the result easily follows.}
\end{proof}
Here, we mention some consequences of the above theorem. Let $M$ be an $s$-subset of $[t]$. For any non-negative integers $i$ and $j$, where $s-w\leq i \leq r$ and $s-r\leq j \leq w$, set
$${\cal F}_i=\{(A^i, B^i) \,\ : \,\  A^i \subseteq M,\  |A^i|=i,\ B^i=M\setminus A^i\},$$
$${\cal E}_j=\{(A^j, B^j) \,\ : \,\  A^j \subseteq M,\  |A^j|=j,\ B^j=M\setminus A^j\}.$$
It is easy to see that $|{\cal F}_i|={s \choose i}$ and $|{\cal E}_j|={s \choose j}$. Also, it is~not difficult to see that ${\cal F}=\cup_{s-w\leq i \leq r}{\cal F}_i$ (resp. ${\cal E}=\cup_{s-r\leq j \leq w}{\cal E}_j$) is a weakly cross-intersecting set-pairs. Therefore, in view of
Theorem~\ref{mainthm}, we have the following corollary which is a generalization of Lemma~\ref{stinlem} (set $s=1$).
\begin{cor}\label{recersive} For any positive integers $0 < s \leq r+w$ and $t\geq r+w$, it holds that
\begin{enumerate}
\item $ N((r,w;d),t)\geq \displaystyle \sum _{s-w\leq i \leq r} {s \choose i}N((r-i,w-s+i;d),t-s)$,
\item $ N((r,w;d),t)\geq \displaystyle \sum _{s-r\leq j \leq w} {s \choose j}N((r-s+j,w-j;d),t-s)$.
\end{enumerate}
\end{cor}
 Let $T((r,w);n)$ denote the maximum number of blocks in an $(r,w)-CFF$ with $n$ points. Erd\"{o}s et al.~\cite{erdos1} discussed $(1,2)$-CFFs in detail, and showed that
$$1.134^n \leq T((1,2);n) \leq 1.25^n.$$
The upper bound is asymptotic and for sufficiently large $n$ is useful. Hence, for large $n$, we have
$N((1,2);t) \geq \frac{1}{\log(1.25)}\log t$.
If we set $s=r+w-3$ in the above corollary, then the following bound can be concluded which can be considered as an improvement of Theorem~\ref{engel}.
\begin{cor}
For any positive integers $r$ and $w$, where $r\geq 2$, it holds that
$$N((r,w),t)\geq {r+w-2 \choose r-1}N((2,1);t-r-w+3)+{r+w-3 \choose r}+{r+w-3 \choose r-3}.$$
\end{cor}

In view of Theorem~\ref{mainthm}, if there exists an $(i,j)$-weakly cross-intersecting set-pairs,
then the following corollary can be concluded. We should mention that
Engel~\cite{engel} obtained a result that is similar to the following corollary.
\begin{cor}\label{engelb} Let $i$, $j$, $r$, and $w$ be positive integers, where $1\leq i\leq r-1$ and $1\leq j\leq w-1$. If there exists an $(i,j)$-weakly cross-intersecting set-pairs of size $g(i,j)$ on a ground set of cardinality $h$, then for any $t$, where $t\geq \max\{h, r+w\}$, we have
$$N((r,w;d),t) \geq g(i,j) N((r-i,w-j;d),t-i-j).$$
\end{cor}
 By a lattice path we mean a path on an $i \times j$ gride from $(0,0)$ to $(i,j)$, where each move is to the right or up. Assume that ${\cal L}(i,j)$ is the set of lattice paths such that the path is strictly below the line $y={j\over i}x$ except at the two endpoints. Tuza \cite{tuza} showed that if $f(i,j)$ is the maximum size of a weakly cross-intersecting set-pairs, then $f(i,j)< \frac{(i+j)^{i+j}}{i^i j^j}$. Recently, Z.~Kir\'{a}ly, Z.L.~Nagy, D.~P\'{a}lv\"{o}lgyi, and M.~Visontai \cite{wif}, by a charming idea and using lattice paths, presented an $(i,j)$-weakly cross-intersecting set-pairs of size $(2i+2j-1)|{\cal L}(i,j)|$ on a ground set of size $2i+2j-1$. Unfortunately, for general $(i,j)$, there is no explicit formula for $|{\cal L}(i,j)|$. However, Bizley \cite{bizley} showed that for relatively prime numbers $i$ and $j$, $|{\cal L}(i,j)|= \frac{{i+j \choose i}}{i+j}$. In \cite{wif}, it is shown $g(i,j) \geq (2-o(1)){i+j\choose i}$, where $f\in o(1)$ means that $\displaystyle\lim_{i+j\rightarrow \infty} f=0$.
\begin{cor} Assume that $r$, $w$, and $t$ are positive integers, where
$t\geq \max\{2r+2w-5, r+w\}$. Then $$N((r,w),t) \geq (2-o(1)) {r+w-2 \choose r-1}{\cal R}(t-r-w+2).$$
\end{cor}
Also, in \cite{wif}, it is shown that there exists an $(r-1,r-1)$-weakly cross-intersecting set-pairs of size $(2-\frac{1}{2r-2}){2r-2 \choose r-1}$ on a ground set of size $4r-6$.
\begin{cor} Assume that $r$ and $t$ are positive integers, where
$t\geq \max\{4r-6, 2r\}$. Then $$N((r,r),t) \geq (2-\frac{1}{2r-2}) {2r-2 \choose r-1}{\cal R}(t-2r+2).$$
\end{cor}
\begin{rem}{\rm It is worth pointing out that the lattice problem is a special case of the generalized ballot problem.
Suppose that in an election, candidate $A$ receives $r$ votes and candidate $B$ receives $w$ votes. Let $r_i$ and $w_i$ denote the number of votes $A$ and $B$ have after counting the $i^{th}$ vote where $1 \leq i \leq r + w$ (notice that $r_i+w_i=i$). Let $k$ be any positive real number. We call a sequence good if $r > kw$ and $r_i > kw_i$ for all $r$. We show the maximum number of good sequence by ${\cal B}(r, w; k)$.  In 1887, Bertrand~\cite{bertrand} showed ${\cal B}(r, w; 1)={r-w\over r+w}{r+w \choose r}$. Determining the exact value of ${\cal B}(r, w; k)$ is known as the generalized ballot problem. It is not difficult to see that ${\cal B}(r,w-1;\frac{r}{w})=|{\cal L}(r,w)|$.
The solution to the generalized ballot problem when $k$ is a positive integer is $ \frac{r-kw}{r+w}{r+w \choose r}$.
Unfortunately, for general $k$, there is no explicit formula for this problem.
In $1962$, Takacs~\cite{takacs} obtained a recurrence formula for the generalized ballot problem. Recently, Delong Meng~\cite{meng} obtained a lower and upper bound for the generalized ballot
problem.}
\end{rem}
The aforementioned bounds improve the existing bounds when the value of $|r-w|$ is small. Now, we present another lower bound
which is an improvement of the earlier bounds whenever $w$ is sufficiently small relative to $r$. Moreover, this bound holds for any $t\geq r+w$.
\begin{thm}\label{estimate} For any positive integers $r$,
$w$, and $t$, where $t\geq r+w$, $r\geq w$, and $r\geq 2$, we have
$$ N((r,w),t) \geq c{{r+w\choose w+1}+{r+w-1 \choose w+1}+ 3 {r+w-4 \choose w-2} \over \log r} \log (t-w+1),$$
where $c$ is a constant satisfies Theorem~\ref{rcff}.
\end{thm}
\begin{proof}{We prove the assertion by induction on $w$.
By Theorem~\ref{rcff}, the assertion holds for $w=1$. Assume that the assertion is
true for every $w'$ where $w'< w$. Easily, one can see that the family
$$\begin{array}{cll}
{\cal F} & = & \{(\emptyset,\{1\}), (\{1\},\{2\}),(\{1,2\},\{3\}), \ldots, \\
         &   & \,\,\, (\{1,2, \ldots, r-w\},\{r-w+1\}), (\{1,2, \ldots, r-w+1\},\{\emptyset\})\}
\end{array}
$$
is a weakly cross-intersecting set-pairs. Hence, in view of Theorem~\ref{mainthm}, it holds that

$$
\begin{array}{ccl}
 N((r,w),t) &  \geq & (\sum_{i=0}^{r-w}N((r-i,w-1),t-i-1))+ N((w-1,w),t-r+w-1)\\
  & = & (\sum_{i=0}^{r-w}N((r-i,w-1),t-i-1))+ N((w,w-1),t-r+w-1)\\
\end{array}
$$
Now by induction we have
$$ \begin{array}{ccl}
N((r,w),t)& \geq & \sum_{i=0}^{r-w}c{{r+w-i-1\choose w}+{r+w-i-2 \choose w}+3{r+w-i-5 \choose w-3}\over \log (r-i)} \log (t-w+1-i)  \\
& & \\
 & + &c{{2w-1\choose w}+{2w-2 \choose w}+3{2w-5 \choose w-3} \over \log (w)} \log (t-r+1).
 \end{array}$$
Since $\frac{\log x}{\log(x-1)}$ is a decreasing function, it holds that
$$
\begin{array}{cll}
 N((r,w),t) & \geq  & c \frac{\log (t-w+1)}{\log r}
 (\sum_{i=0}^{r-w}{r+w-i-1\choose w}+{r+w-i-2 \choose w}+3{r+w-i-5 \choose w-3}) \\
& & \\
& + & c \frac{\log (t-w+1)}{\log r}({2w-1\choose w}+{2w-2 \choose w}+3{2w-5 \choose w-3})\\
& & \\
 & \geq  & c \frac{\log (t-w+1)}{\log r}
 (\sum_{i=0}^{r-w}{r+w-i-1\choose w}+{r+w-i-2 \choose w}+3{r+w-i-5 \choose w-3}) \\
& & \\
& + & c \frac{\log (t-w+1)}{\log r}({2w-1\choose w+1}+{2w-2 \choose w+1}+3{2w-5 \choose w-2})\\
& & \\

& = & c{{r+w\choose w+1}+{r+w-1 \choose w+1}+3{r+w-4 \choose w-2}\over \log r} \log (t-w+1).
\end{array}
$$
}
\end{proof}
\section{Fractional Biclique Cover} The next result concerns the
fractional version of biclique cover. If ${\cal R}$ is the set of
all bicliques of a graph $G$, then each biclique cover of $G$ can
be described by a function $\phi:{\cal R}\rightarrow\{0,1\}$ such
that $\phi(G_i)=1$ if and only if $G_i$ belongs to the cover.
Hence, $bc(G)$ is the minimum of $\displaystyle \sum_{G_i \in
{\cal R}}\phi(G_i)$ over all function $\phi:{\cal
R}\rightarrow\{0,1\}$ such that for any edge $e$ of $G$ we have
\begin{equation}
\label{fraction} \displaystyle \sum_{G_i \in {\cal R}: e \in
E(G_i)}\phi(G_i)\geq 1.
\end{equation}
The fractional biclique covering number $bc^*(G)$ is the minimum
of $\displaystyle \sum_{G_i \in {\cal R}}\phi(G_i)$ over all
functions $\phi:{\cal R}\rightarrow[0,1]$  satisfying
(\ref{fraction}).\\
Fractional graph theory is the modification  of integer-valued
graph parameters  to take its value on non-integer values. For
more on fractional graph theory and other fractional graph
parameters, see \cite{fractional}.
 In the fractional cover, using linear programming, it is proved that
 $$bc^*(G)=\displaystyle \inf_d \frac{bc_d(G)}{d}=\displaystyle\lim_{d \rightarrow\infty}
 \frac{bc_d(G)}{d}.$$
 Also, we have the following theorem.
\begin{alphthm}\label{frac}{\rm \cite{fractional}}
 For every non-empty edge-transitive graph $G$, we have
$$bc^*(G)=\frac{|E(G)|}{B(G)},
$$
where $B(G)$ is the maximum  number of edges among the bicliques
of $G$.
\end{alphthm}
 Easily, one can see that
$$B(I_t(r,w))= \displaystyle \max_{t'+t''=t} {t'\choose r}
{t''\choose w}.$$

Also, we have $|E(I_t(r,w))|={t\choose r}{t-r \choose w}$, and
 $I_t(r,w)$ is an edge-transitive graph. Therefore, in view of Theorem~\ref{frac}, we have
\begin{center}
$bc^*(I_t(r,w))=\displaystyle\min_{t'+t''=t}{{t\choose r}
{t-r\choose w} \over {t'\choose r}{t''\choose w}}$.
\end{center}
By a straightforward calculation, one can see that

$$bc^*(I_t(r,w))=\displaystyle\min_{t'+t''=t}{{t\choose r} {t-r\choose w} \over
{t'\choose r}{t''\choose w}} = \displaystyle \min_{w\leq m \leq
t-r} {{t\choose m}\over {t-r-w \choose m-w}}.$$

Lov\'{a}sz \cite{lovasz} proved that for any graph $G$ with
maximum degree $\Delta(G)$
 $$bc^*(G)\geq\frac{bc(G)}{1+\ln(\Delta(G))}.$$

The maximum degree of  the graph $I_t(r,w)$ is equal to
 $$\max\{{t-w\choose r},{t-r\choose w}\}={t-w \choose r}.$$ So we have the following
 corollary.
 \begin{cor}For any  positive integers $r$, $w$, and $t$, where $t\geq
 r+w$, we have
 $$N((r,w),t)\leq \displaystyle\min_{w\leq m \leq
t-r} {{t\choose m}\over {t-r-w \choose m-w}}(1+\ln({t-w \choose
r})).$$
 \end{cor}
In \cite{engel}, Engel  proved that
 $$ N((r,w),t)\geq \displaystyle\min_{w-1\leq m \leq
t-r+1} {{t\choose m}\over {t-r-w+2 \choose
m-w+1}}(N((1,1),t-r-w+2)).$$
 Hence, we have
$$ N((r,w),t)\geq \displaystyle\min_{w-1\leq m \leq
t-r+1} {{t\choose m}\over {t-r-w+2 \choose m-w+1}}(\log_2
(t-r-w+2) +\frac{1}{2}\log_2\log (t-r-w+2) +c),$$

where $c$ is a constant. In the next theorem, we specify the exact
value of $N((r,w;d),t)$ for some special value of $d$. In the
proof of the next theorem, by $S_t$ we mean the permutation group
of the set $[t]$.
\begin{thm}\label{equality}
For any positive integers $r$, $w$, $t$, $d_0$, and $d =
\frac{B(I_t(r,w))}{|E(I_t(r,w))|} t!$, where
$t\geq r+w$, we have
$$N((r,w;d_0d),t)=d_0 (t!).$$
\end{thm}
\begin{proof}{For every
$\sigma \in S_t$, define the function $ f_\sigma :
V(I_t(r,w))\rightarrow V(I_t(r,w)) $ such that for every set $A =
\{ i_1, i_2, \ldots, i_l \}\in V(I_t(r,w))$, we have
$f_\sigma(A)=\{\sigma(i_1), \ldots, \sigma(i_l) \}$ (note that
here $l=r$ or l=w). Set $G=\{f_\sigma \,\ | \,\,\ \sigma \in S_t
\}$. One can see that $G$ is a subgroup of $Aut( I_t(r,w) )$ and
also $G$ acts transitively on $E(I_t(r,w))$. Now it is simple to
check that
$$ {bc_d(I_t(r,w))\over d} = \frac{|E(I_t(r,w))|}{B(I_t(r,w))}.$$
To see this, assume that $K$ is a biclique of $I_t(r,w)$, where
$|E(K)|=B(I_t(r,w))$. Construct a biclique cover of $I_t(r,w)$ as
follows. Set
$$ {\cal C}= \{f_\sigma(K) \,\,\ | \,\,\ \sigma \in S_t\}.$$
It is readily seen that ${\cal C}$ is a biclique cover and
every edge is covered with exactly
$d=\frac{B(I_t(r,w))t!}{|E(I_t(r,w))|}$ bicliques. So
 $${bc_d(I_t(r,w))\over d}\leq {|E(I_t(r,w))| \over B(I_t(r,w))}.$$
 On the other hand, by the definition of fractional biclique cover, for every graph $G$ and
  every positive integer $d$
 we have $bc^*(G)\leq {bc_d(G)\over d}$. Particularly,
 $$bc^*(I_t(r,w))\leq {bc_d(I_t(r,w))\over d}.$$
 Consequently, in view of Theorem~\ref{frac}, we have

 $${bc_d(I_t(r,w))\over d} = {|E(I_t(r,w))| \over B(I_t(r,w))}.$$

Also, for any positive integer $d_0$ we have
$$ bc_{d_0d}(I_t(r,w)) \leq d_0 bc_{d}(I_t(r,w)).$$
Hence,
\begin{equation}
 \label{equal}
 {|E(I_t(r,w))| \over B(I_t(r,w))} = bc^*(I_t(r,w))\leq \frac{bc_{d_0d}(I_t(r,w))}{d_0d} \leq
 \frac{bc_{d}(I_t(r,w))}{d} = {|E(I_t(r,w))| \over B(I_t(r,w))}.
\end{equation}
Consequently, using (\ref{equal}) we obtain the result. }
\end{proof}
\begin{cor} For any positive integers $r$, $w$, and $c$, we have
$$N((r,w;c(r+1)!w!), r+w+1)=c(r+w+1)!.$$
\end{cor}
\begin{proof}{Set $G:=I_{r+w+1}(r,w)$. The
 graph $G$ is $C_4$ free. So $B(G)=r+1$.
Also, $|E(G)|= (w+1){r+w+1 \choose r}$. Hence, by a straightforward calculation
and using Theorem~\ref{equality}, the corollary follows.  }
\end{proof}
An $n \times n$ matrix $H$ with entries $+1$ and $-1$ is called a {\it Hadamard matrix} of {\it order} $n$
if $HH^t = nI$. It is seen that any two distinct columns of $H$ are orthogonal.
Also, if we multiply some rows or columns by $-1$, or if we permute rows or columns, then $H$ is still
a Hadamard matrix. Two such Hadamard matrices are
called equivalent. Easily, for any Hadamard matrix $H$, we can find an
equivalent one for which the first row and the first column consist
entirely of $+1$'s. Such a Hadamard matrix is called {\it normalized}.
\begin{thm}Let $d$ be a positive integer such that there exists a Hadamard matrix of order $4d$, then
$N((1,1;d),4d-1)=4d-1.$
\end{thm}
\begin{proof}{Let $H$ be a normalized Hadamard matrix of order
$4d$. Delete the first row and the first column. Also, assume that $K^-_{4d-1,4d-1}$ has $(X,Y)$ as its vertex set where $X=\{v_1, \ldots, v_{4d-1}\}$ and $Y=\{v'_1, \ldots, v'_{4d-1}\}$.
 Assign to the $j$th column of $H$, two sets $X_j$ and $Y_j$ as follow
$$ X_j=\{ v_i | h_{ij}=+1 \}\quad \&  \quad Y_j=\{v'_i | h_{ij}=-1 \}.$$
Construct a complete bipartite graph $G_j$ with vertex set $(X_j, Y_j)$. The edge $v_iv'_j$ is covered by the complete bipartite graph $G_k$ if and only if the corresponding entries of column $k$ in row $i$ is $+1$ and in row $j$ is $-1$. It is well-known that the number of these columns, in a normalized Hadamard matrix of order $4d$, is equal to $d$. Hence, every edge is covered exactly $d$ times. So $bc_{d}(K^-_{4d-1,4d-1})\leq 4d-1$. On the other hand, $K^-_{4d-1,4d-1}$ is an edge-transitive graph. Therefore, in view of Theorem~\ref{frac}, we have $${4d-1\over d}=\frac{|E(K^-_{4d-1,4d-1})|}{B(K^-_{4d-1,4d-1})}\leq {bc_d(K^-_{4d-1,4d-1})\over d}.$$ Consequently, $ 4d-1 \leq bc_{d}(K^-_{4d-1,4d-1})$ and
the result follows.
}
\end{proof}
\def\cprime{$'$}

\end{document}